\documentclass[12pt]{article}
\usepackage{verbatim}
\usepackage{amsfonts}
\usepackage{graphics}
\usepackage{amsmath}
\usepackage{times}
\usepackage{appendix}
\usepackage{color}
\usepackage{enumerate}
\usepackage{fancyhdr,latexsym,amsmath,amsfonts,amssymb,amsbsy,amsthm,url}
\usepackage[margin=0.5in,footskip=0.25in]{geometry}
\usepackage{graphics,graphicx,epsfig}
\usepackage{breqn}
\usepackage{caption,subcaption,float,subfloat}
\usepackage{pdflscape}
\usepackage{hyperref}
\usepackage[ruled,vlined]{algorithm2e}
\newtheorem{theorem}{Theorem}[]

\usepackage[english]{babel}
\usepackage[dvipsnames]{xcolor}
\usepackage{tikz}

\numberwithin{equation}{section}

\fancypagestyle{usual}{
	\lhead[\thepage]{}
	\chead[{\it \shortauthors}]{\sc \shorttitle}
	\rhead[]{\thepage}
	\lfoot[]{}
	\cfoot[]{}
	\rfoot[]{}
	
}

\makeatletter

\renewcommand{\section}{
	\@startsection
	{section}% name
	{1}% level
	{0pt}% indent
	{1.1\baselineskip}% beforeskip
	{0.2\baselineskip}% afterskip
	{\sc \centering}% style
}

\renewcommand{\subsection}{
	\@startsection
	{subsection}% name
	{1}% level
	{0pt}% indent
	{1.1\baselineskip}% beforeskip
	{0.2\baselineskip}% afterskip
	{\sc \centering}% style
}

\renewcommand{\subsubsection}{
	\@startsection
	{subsubsection}% name
	{1}% level
	{0pt}% indent
	{1.1\baselineskip}% beforeskip
	{0.2\baselineskip}% afterskip
	{\sc \centering}% style
}

\makeatother

\usepackage[flushleft]{threeparttable}
\usepackage{rotating,booktabs,multirow}
\usepackage{ colortbl}
\usepackage{makecell, cellspace, caption}

\begin{document}
	\everymath{\displaystyle}
		\title{\large\sc Continuous-Time Higher Order Markov Chains : Formulation and Parameter Estimation}
	\normalsize
	\author{\sc{Suryadeepto Nag} \thanks{Indian Institute of Science Education and Research Pune, Pune-411008, Maharashtra, India, e-mail: suryadeepto.nag@students.iiserpune.ac.in}
		}
	\date{}
	\maketitle
\begin{abstract}
	Stochastic processes find applications in modelling systems in a variety of disciplines. A large number of stochastic models considered are Markovian in nature. It is often observed that higher order Markov processes can model the data better. However most higher order Markov models are discrete. Here, we propose a novel continuous-time formulation of higher order Markov processes, as stochastic differential equations, and propose a method of parameter estimation by maximum likelihood methods. 
\end{abstract}
\section{Introduction}
While Markov chains have been used in a variety of disciplines such as biology and economics \cite{Norris97}, the usage of higher order Markov chains is relatively sparse.  This may be attributed to the significantly larger number of parameters one encounters in higher order Markov models \cite{Raftery85}. However, the Markovian nature of many models often fail to stand the test of data. In \cite{Marathe05}, Marathe and Ryan examined the validity of modelling sets of data as Markovian geometric brownian motions, and they observed instances where there exists a dependence of increments in the state variable on the past increments. In \cite{Shorrocks76}, Shorrocks investigated the Markovian assumption in modelling income mobility and concluded that a second order model would be more appropriate.  \\
\\
 Another disadvantage of using a higher order Markov chain in modelling so far, has been in the discrete nature of the models. We observe that most instances of higher order Markov models which have been used so far, have involved discrete time models. This leaves out a very large class of stochastic processes which cannot be modelled using higher order Markov models as they are inherently continuous in nature. Discretizing such systems not only diminishes the accuracy of modelling, but also significantly increases the ``order" of the now discrete model, in turn increasing the number of parameters in the model. Despite these factors, higher order Markov chains have elicited sustained interest \cite{Ching08}. \\\\
The objective of this paper is to develop a continuous-time formulation of higher order Markov chains. The inclusion of continuous processes under higher order Markov models will make it possible to model better, several systems, which were till now modelled either by Markov processes or by discrete-time higher order Markov processes.  Under our formulation, continuous-time higher order Markov processes are modelled with the order being an interval of time instead of a number of discrete states. The state variable is written as a function of an underlying time variable, such that the ``state function" represents the state variable over all instances of time within the interval with length equal to the order. Drift and diffusion functions are also written in terms of the underlying variable. This allows higher order Markov processes to be modelled as stochastic differential equations, as is done with Markov processes. We subsequently derive the Fokker-Planck equation in terms of these state functions which lets us find the transition probability from one state function to another. This can be used to estimate parameters making use of a maximum-likelihood approach. 

\section{A novel formulation of Continuous time higher order Markov processes (chains)} 
We know that a higher order Markov process (HOMP) is a stochastic process that depends on not just the current state but also the previous states. Now consider a discrete time HOMP of order $n=\frac{\tau}{\Delta t}$ that depends on all states of the system in the last $\tau$ time, such that the time $\tau$ is divided into equal intervals of length $\Delta t$ (one step). \\
Accordingly, we can say that 
\begin{equation}
x(T+\Delta T)=f(x(T),x(T-\Delta T),\dots,x(T-\tau))
\end{equation}
This becomes a continuous time Markov process at the $0$ limit of $\Delta t$, i.e., 
\begin{equation}
\lim\limits_{\Delta T\to 0}x(T+\Delta T)=f(x(T),x(T-\Delta T),\dots,x(T-\tau))
\end{equation}
Now let us define $F$ as some function over the $t-$space, defined over the interval $[T-\tau,T]$, such that for some underlying time variable $x$, the function $F_t(x)$ defines a HOMP of order $\tau$, with $x$ taking values in $[t-\tau,t]$.  Defining the HOMP in this fashion allows for us to approach the problem analogously to a Markov process, where every new state is a function of only the state variable. It may be noted that the state function, as defined here would not only be the state at the present time, but also include all states in the last $\tau$ time. This formulation will help derive the Fokker-Planck equation which will subsequently prove useful in the parameter estimation. Thus for some transformation $V$ we can define a HOMP as, 
\begin{equation*}
\lim\limits_{\Delta t\to 0}F_{t+\Delta t}(x)=V(F_t(x))
\end{equation*}
Let $h(t)$ be some function such that $h(t)\to 0\ \forall\ t\in[T,T-\tau]$. Therefore we can define a derivative of the transformation, 
\begin{equation*}
V'(F,h)=\frac{V(F+h)-V(F)}{h}
\end{equation*}
We can now generalise this definiton of $V'$ for any $h(x)$ regardless of its form. Consider the function $V(F)=F^a$
\begin{align*}
\therefore V'(F,h)&=\lim\limits_{h\to 0}\frac{(F+h)^a-(F)^a}{h}\\
&=\lim\limits_{h\to 0}\frac{F^a+aF^{a-1}h+...+h^a-F^a}{h}\\
&=\lim\limits_{h\to 0}\frac{aF^{a-1}h+{a\choose 2}F^{a-2}h^2+...+h^a}{h}\\
&=\lim\limits_{h\to 0}aF^{a-1}+h(G(a,h))\\
&=aF^{a-1}.
\end{align*}
Now, since the derivative of a transformation of the form $F^a$ doesn't depend on the form of $h$, any transformation can be written as an infinite power series of $F^a$. Also, it is shown that $V'(F,h)$ does not depend on the form of $h$. Henceforth, we shall denote the derivative of a transformation $V$ with respect to $F$ by $V'(F)$.
\\\\
We can define the Integral $\int_{F} V(F)\ dF$ of $V(F)$ with respect to $F$ in Riemannian terms by integrating over the $F-$space.
\section{Stochastic differential equations}
\subsection{Formulation}
We can write a HOMP of order $\tau$ as a stochastic differential equation (SDE), 
\begin{equation}
dy(t)=\nu(H_t)dt+\varsigma(H_t)dW,
\end{equation}
where $y(t)$ is the state variable, and $H_t$ is the state  function that represents the past realisations of the state variable from $t-\tau$ to $t$, with $dW$ being a standard Wiener process. Note that here, $\nu$ and $\varsigma$ are scalar quantities. We would like to rewrite this equation such that $H_t$ is the state function. For this, we define an underlying variable $x\in \mathbb{R}$, and we define the functions, $\mu(H_t(x))$, $\sigma(H_t(x))$ and $dB(x)$, along with $H_t(x)$ in the time dependent domain $x\in [t-\tau,t]$ for all of them. Here $\mu$ is the drift function, $\sigma$ is the diffusion function and $dB(x)$ is the $x$ dependent Wiener process. Therefor,e we can rewrite the SDE as,
\begin{equation}
dH_t(x)=\mu(H_t(x))dt+\sigma(H_t(x))dB(x)
\end{equation}
While the function $dB$ could have been written as $dB(H_t(x))$, it has been defined only in $x$ and independent of $H_t$, because $dB(H_{t_1}(x_0))=dB(H_{t_2}(x_0))$. Eq.(3.1) is related to eq.(3.2) as the realizations of the functions on either side of eq. (3.2) at $x=t$ gives eq. (3.1).
Note that unlike eq. (3.1), each term on the right hand side of eq. (3.2) returns not a scalar, but a function in $x$. The state function therefore contains the entire information in itself, thereby, in principle, reducing the HOMP process to a Markov process. 
\subsection{Examples}
Here are some simple examples of continuous time higher order Markov processes 
\begin{enumerate}
	\item \textbf{A higher order Geometric Brownian Motion}\\
	\begin{equation}
	dy(t)=\alpha \int_{t-\tau}^{t}H_t(x)dx\ dt + \beta \int_{t-\tau}^{t}H_t(x)dx\ dW(t),
	\end{equation} 
	where $y(t)=H_t(t)$.
	This system is analogous to the Markovian brownian motion SDE
	\begin{equation*}
	dy(t)=\alpha y\ dt +\beta y\ dW(t)
	\end{equation*}
	In the formulation developed in this paper, this can be written as
	\begin{equation}
	dH_t(x)=\alpha \int_{x-\tau}^{x}H_t(x')dx'\ dt+\beta \int_{x-\tau}^{x}H_t(x')dx'\ dB(x)
	\end{equation}
	Here $\mu(H_t(x))= \alpha \int_{x-\tau}^{x}H_t(x')dx'$ is the drift term and $\sigma=\beta \int_{x-\tau}^{x}H_t(x')dx'$ is the diffusion term. Observe that fixing $x=t$ in eq. (3.3) gives us eq.(3.4). This is a simple higher order Markov process with only 2 parameters. 
	\item \textbf{A higher order Orstein-Uhlenbeck process}
	\begin{equation}
	dy(t)=-\theta \int_{t-\tau}^{t}H_t(x)dx\ dt +\sigma dW
	\end{equation}
	or 
	\begin{equation}
		dH_t(x)=-\theta \int_{x-\tau}^{x}H_t(x')dx'\ dt +\sigma dB(x)
	\end{equation}
	This can be used to model systems with mean-reversions, where the reversion of the mean depends on the moving average of the state variable, instead of just the present value of the state variable. This helps explain mean reversions with greater time periods of reversion. \\\\
	Both the examples discussed till now depend on moving averages of the state variable, rather than simply the present state. This may open doors to new avenues of stochastic modelling, especially in finance, where moving averages are of particular interest. This may find particular application in commodity pricing, along the lines of models explored by Wets and Riot \cite{Wets15} and Nag et al \cite{Nag20}.
	\item \textbf{Weighted Moving Average Models}\\
	The examples so far have explored models with a finite (and rather small) number of parameters. However, often more realistic models would depend on a much larger number of variables. For instance, Hull \cite{Hull18}, discusses modelling volatility from moving averages. A stochastic differential model based on this would be of the form, 
	\begin{equation}
	dy(t)=\nu(H_t)dt + \bigg(\int_{t-\tau}^{t}\omega(x)\sigma^2(x)dx\bigg)^{\frac{1}{2}}\ dW(t)
	\end{equation}
	or 
	\begin{equation}
	dH_t(x)=\mu(H_t(x))dt + \bigg(\int_{x-\tau}^{x}\omega(x')\sigma^2(x')dx'\bigg)^{\frac{1}{2}}\ dB(x)
	\end{equation}
	Here the diffusion term is written as a function of historic values of the itself i.e. 
	\begin{equation*}
	\sigma^2(x)=\int_{x-\tau}^{x}\omega(x')\sigma^2(x')dx'
	\end{equation*}
	Further, $\omega(x)$ is the parameter function, physically representing weights given to the different historical values of volatility. This is useful in predictive modelling as it can be used to weigh more recent values of volatility higher. In particular, a continuous time model based on the Exponentially Weighted Moving Average, used to model volatility \cite{Kuen92,Ederington05} can be written as, 
	\begin{equation*}
	\sigma^2(x)=\int_{x-\tau}^{x}\lambda^{x'-(t-\tau)}\sigma^2(x')dx'.
	\end{equation*}	
	Here, the parameter function $w_t(x)$ at any time $t$ can be written as, 
	\begin{equation*}
	w_t(x)=\lambda^{x'-(t-\tau)}.
	\end{equation*}
\end{enumerate}
\section{The Chapman-Kolmogorov and the Fokker-Planck equations} 
In this section,  we derive the Chapman-Kolmogorov equation, and subsequently derive the Fokker-Planck equation under our framework. This allows us to get an equation that describes the evolution of the transition probability, solving which gives the transition probability of transitioning from one state function to another, which can later be used in parameter estimation. 
\subsection{The Chapman Kolmogorov equation} 
\begin{theorem}
	For a $\tau^{th}$ order Markov process,
	\begin{equation}
	P\{G_T|F_0\}=\int_{H_t} P\{G_T|H_t\}P\{H_t|F_0\}dH_t
	\end{equation}
	\\
	where $G_T$ and $H_t$ are the curves that denote the state of the system over time $(T-\tau,T)$ and $(t-\tau,t)$ respectively and $F_0$ is defined over $(-\tau,0)$ such that $0<t<T$ and $P\{B|A\}$ is the conditional probability of event $B$ occuring given $A$ has occured.
\end{theorem}
\begin{proof}
	Here,
	\begin{align*}
	&P\{G_T|F_0\}\\
	=&P\{G_T=g|F_0=f\}\\
	=&\int_{H_t} P\{(G_T=g)\cap (H_t=h)|F_0=f\}dH_t\\
	=&\int_{H_t} P\{G_T=g|(H_t=h)\cap(F_0=f)\}P\{H_t=h|F_0=f\} dH_t\\
	=&\int_{H_t} P\{G_T=g|H_t=h\}P\{H_t=h|F_0=f\}dH_t\\
	=&\int_{H_t} P\{G_T|H_t\}P\{H_t|F_0\}dH_t
	\end{align*}
\end{proof}

Assuming uniform order and time homogeneity, we have written the notation in a manner such that each state function is only denoted by the upper interval of the time it is defined over in the subscript. 
\subsection{The Fokker-Planck equation}
We now derive the Fokker-Planck Equation for HOMPs analogous to the proof in Kolpas \cite{Kolpas08} and Adesina \cite{Adesina08} for Markov processes, the only difference being our states are described by functions (rather than variables) with transformations acting on them. \\\\
\begin{theorem}
	\begin{equation}
	\frac{\partial P(H_t|F_0)}{\partial t}=-\frac{\partial}{\partial H_t}\big[V(H_t)P(H_t|F_0)\big]+\frac{\partial^2}{\partial H_t^2}\big[D(H_t)P(H_t|F_0)\big],
	\end{equation}
	where $V(H_t)=D^{(1)}(H_t)$ and $D(H_t)=D^{(2)}(H_t)$, with 
	\begin{equation*}
	D^{(n)}(H_t)=\frac{1}{n!}\lim\limits_{\Delta t\rightarrow0}\frac{1}{\Delta t}\int_{G_t} (G_t-H_t)^nP(G_{t+\Delta t}|H_t)dG_t,
	\end{equation*}
	where $H_t$ is a higher order Markov Process of order $\tau$.
\end{theorem}
\begin{proof}
	Let us consider a smooth transformation with compact support $\phi$. Now, consider the integral,
	\begin{equation*}
	\int_{G_t}\phi(G_t)\frac{\partial P(G_t|F_0)}{\partial t}dG_t,
	\end{equation*}
	where $F_0=F(0)$ is a curve denoting the state of the process over the time $[-\tau,0]$, and $G_t=G(t)$ is a curve denoting the state of the process over the time $[t-\tau,t]$. We shall use this notation to denote all functions defined over a time of length $\tau$. \\From the definition of partial derivatives, we know that, 
	\begin{equation*}
	\frac{\partial P(G_t|F_0)}{\partial t}=\lim\limits_{\Delta t\to 0}\frac{P(G_{t+\Delta t}|F_0)-P(G_t|F_0)}{\Delta t}.
	\end{equation*}
	Therefore we have, 
	\begin{equation*}
	\int_{G_t}\phi(G_t)\frac{\partial P(G_t|F_0)}{\partial t}dG_t=\lim\limits_{\Delta t\to 0}\int_{G_t}\phi(G_t)\frac{P(G_{t+\Delta t}|F_0)-P(G_t|F_0)}{\Delta t}dG_t
	\end{equation*}
	Applying the Chapman-Kolmogorov equation as derived (in the preceding subsection), to the right hand side, we get,

	\begin{equation*}
	\lim\limits_{\Delta t\rightarrow0}\frac{1}{\Delta t}\bigg[\int_{G_t}\phi(G_t)\int_{H_t}P(G_{t+\Delta t}|H_t)P(H_t|F_0)dH_tdG_t-\int_{G_t}\phi (G_t)P(G_t|F_0)dy\bigg]
	\end{equation*}

	Using \ \ \ $\int_{G_t} P(G_{t+\delta t}|H_t)dG_t=1$, we obtain the expression as 
	\begin{align*}
	\lim\limits_{\Delta t\rightarrow0}\frac{1}{\Delta t}\bigg[&\int_{H_t} P(H_t|F_0)\int_{G_t} P(G_{t+
		\Delta t}|H_t) \phi(G_t)dG_tdH_t\\-&\int_{G_t} P(G_{t+\Delta t}|H_t)\int_{G_t} \phi(G_t) P(G_t|F_0)dG_tdG_t\bigg]
	\end{align*}

	Replacing $G_t$ with $H_t$ in the second term (since it is integrated over all curves anyway), we get,
	\begin{equation*}
	\lim\limits_{\Delta t\rightarrow0}\frac{1}{\Delta t}\bigg[\int_{H_t} P(H_t|F_0)\int_{G_t} P(G_{t+
		\Delta t}|H_t)(\phi(G_t)-\phi(H_t))dG_tdH_t\bigg]
	\end{equation*}

	Taking a Taylor expansion of $\phi(G_t)$ about $H_t$ gives,
	\begin{equation*}
	\lim\limits_{\Delta t\rightarrow0}\frac{1}{\Delta t}\bigg[\int_{H_t} P(H-t|F_0)\int_{G_t} P(G_{t+
		\Delta t}|H_t)\sum_{n=1}^{\infty}\phi^{(n)}(H_t)\frac{(G_t-H_t)^n}{n!}dG_tdH_t\bigg].
	\end{equation*}
	Defining the jump moments as,
	\begin{equation*}
	D^{(n)}(H_t)=\frac{1}{n!}\lim\limits_{\Delta t\rightarrow0}\frac{1}{\Delta t}\int_{G_t} (G_t-H_t)^nP(G_{t+\Delta t}|H_t)dG_t,
	\end{equation*}
	it follows that,
	\begin{equation*}
	\int_{G_t} \phi(G_t)\frac{\partial P(G_t|F_0)}{\partial t}dG_t=\int_{H_t} P(H_t|F_0)\sum_{n=1}^{\infty}D^{(n)}(H_t)\phi^{(n)}(H_t)dH_t
	\end{equation*}
	Integrating the right side by parts $n$ times and shifting to the left,
	\begin{equation*}
	\int_{H_t} \phi(H_t)\biggl(\frac{\partial P(H_t|F_0)}{\partial t}-\sum_{n=1}^\infty\bigg(-\frac{\partial}{\partial H_t}\bigg)^n\big[D^{(n)}(H_t)P(H_t|F_0)\big]\biggr)dH_t=0.
	\end{equation*}
	Since $\phi(X)$ is arbitrary, it follows that,
	\begin{equation*}
	\frac{\partial P(Ht_|F_0)}{\partial t}-\sum_{n=1}^\infty\bigg(-\frac{\partial}{\partial H_t}\bigg)^n\big[D^{(n)}(H_t)P(H_t|F_0)\big]=0.
	\end{equation*}
	Setting $D^{(n)}(H_t)=0$ for $n>2$, we get,
	\begin{equation*}
	\frac{\partial P(H_t|F_0)}{\partial t}=-\frac{\partial}{\partial H_t}\big[V(H_t)P(H_t|F_0)\big]+\frac{\partial^2}{\partial H_t^2}\big[D(H_t)P(H_t|F_0)\big],
	\end{equation*}
	where $V(H_t)=D^{(1)}(H_t)$ and $D(H_t)=D^{(2)}(H_t)$.
\end{proof}

This is the Fokker-Plack Equation for the HOMP.
From the definition of $D^{(1)}(X)$ and $D^{(2)}(X)$ it follows that they are the drift and the diffusion terms respectively for the HOMP.
\subsection{Smooth processes}
For a smooth process, i.e. where $H_t$ is a smooth function of $x$, we can expand $H_t$ into its Taylor expansion.
\begin{equation*}
H_t(x)=c_0 +c_1 x+c_2x^2+\dots
\end{equation*}
Now, we can further write functions of $H_t$ as functions of $x$, and subsequently simplify the Fokker-Planck equation. We make the following substitutions, 
\begin{equation*}
P(H_t(x)|F_0)=p(x) , \ \ \ V(H_t)=v(x), \ \ \, D(H_t)=\varphi(x).
\end{equation*}
Therefore,
\begin{equation*}
 \partial_ {H_t}\big[V(H_t)P(H_t|F_0)\big]=\frac{\partial_x[v(x)p(x)]}{\partial_x H_t(x)}
\end{equation*}
and
\begin{equation*}
\partial^2_{ H_t}\big[D(H_t)P(H_t|F_0)\big]= \bigg(\partial_xH_t(x)\bigg)^{-1}.\partial_x\biggl[\frac{\partial_x\varphi(x)p(x)]}{\partial_xH_t(x)} \biggr]
\end{equation*}
Thus we get the Fokker-Planck equation in $x$ as,
\begin{equation}
\partial_tp(x,t)=\bigg(\partial_xH_t(x)\bigg)^{-1}.\biggl[\partial_x[v(x)p(x,t)]+\partial_x\bigg(\frac{\partial_x[\varphi(x)p(x,t)]}{\partial_xH_t(x)} \bigg) \biggr]
\end{equation}
The Fokker-Planck equation is the continuity equation for probability. Therefore, it can also be written as, 
\begin{equation*}
\partial_tp=\partial_{H_t}j
\end{equation*}
At equilibrium, we know that the probability current must be 0. Therefore, we get the following condition, 
\begin{equation*}
v(x)p(x,t)+\frac{\partial_x[\varphi(x)p(x,t)]}{\partial_xH_t(x)}=0.
\end{equation*}
Furthermore, analogising from the Fokker-Planck equation for Markov processes \cite{Peters17}, we get, 
\begin{equation*}
v(x)=\mu(H_t(x)) ~\text{and}~ \varphi(x)=\sigma(H_t(x))^2
\end{equation*}
\section{Parameter estimation} 
As is the case with most problems in mathematical modelling, the estimation of the parameters are of particular interest in ascertaining the validity of the model and in also gaining insights about physical quantities that are represented by or related to the parameters.
Maximum likelihood estimation in Markov models and SDE models have been investigated before \cite{Kleinhans07,Lindstrom07,Kaltenbacher18}. Similar approaches of using the transition probability to find the likelihood can be extended to HOMPs.
Parameter estimation in classical model of a discrete time HOMP of order $n$ usually depends on $(m-1)m^n$ parameters, where $m$ is the number of possible states \cite{Ching03}. The estimation of such a large number of parameters involves large matrices and is therefore computationally heavy. This problem is further magnified in the case of continuous processes, as any reasonably good discretisation would make the order of the then discrete time process very large, and the estimation - computationally expensive.
The best possible approach, designed to circumvent this problem, is by approximating the parameter function, by a fewer number of parameters, and approximating the functions for drift and diffusion. We propose a method for the estimation of parameters for a specific class of problems, where the drift and diffusion functions $\mu$ and $\sigma$ are smooth. While, strictly, there may be a large number of systems that do not belong to this set, the set will still consist of most real systems, as we would expect the function to act similarly to points in the neighbourhood i.e., 
\begin{equation*}
\lim\limits_{x\to x'}\mu(H_t(x))=\mu(H_t(x')), \ \ \ \ \ \lim\limits_{x\to x'}\mu'(H_t(x))=\mu'(H_t(x'))
\end{equation*}
\begin{equation*}
\lim\limits_{x\to x'}\sigma(H_t(x))=\sigma(H_t(x')), \ \ \ \ \ \lim\limits_{x\to x'}\sigma'(H_t(x))=\sigma'(H_t(x'))
\end{equation*}
The advantage of considering smooth or piecewise-smooth functions is that they allow us to use polynomials or splines to approximate the parameter functions within the drift/diffusion terms in a significantly smaller number of constant parameters. 
Additionally, while we may believe a continuous model may be the best explanation of a system, continuous data is almost unheard of. Most available data is discrete, and at times irregular (for instance stock market data has irregular numbers of working days in a week). The discrete data can be interpolated and made continuous for the purpose of estimating model parameters. This also helps deal with irregularities since, if the number of parameters used to interpolate the data are sufficiently small, irregularities in the availability of data will be almost inconsequential. 
Once, the data $Y$ has been partitioned into $n$ intervals (in time) of length $\tau$ (The last interval can be shorter) and each piece $H_t$ interpolated in terms of the parameters $\Theta_d(H_t)$, and we have a guess for the approximated model parameters $\Theta_m$, we can calculate the transition probability after solving the Fokker-Planck equation. Therefore for $t_1,t_2,t_3,...,t_n$ defined as $\tau, 2\tau,3\tau,...,n\tau$ we can calculate the likelihood function by, 
\begin{equation*}
\mathcal{L}(Y,\Theta_m)=p(H_{t_n}|H_{t_{n-1}}).p(H_{t_{n-1}}|H_{t_{n-2}})....p(H_{t_3}|H_{t_2}).p(H_{t_2}|H_{t_1}).
\end{equation*}
Consequently, we can use this likelihood function to arrive at the Maximum-likelihood estimate of $\Theta_m$. 
\section{Conclusion}
In this paper, we have have introduced a new formulation for the treatment of continuous-time higher order Markov processes as stochastic differential equations.  The advantage of such a formulation allows better modelling of a large number of continuous processes which either had to be modelled as first order Markov processes or as discrete processes. The Chapman-Komogorov and the Fokker-Planck equations are derived and, a method for estimating parameters is laid out. The problem of a large number of parameters typically observed in higher order Markov process is addressed by alluding to the likely smooth nature of parameter functions, which allows them to be approximated reasonable by a sufficiently small number of  constant parameters.
\section{Acknowledgements}
We would like to thank Professor Siddhartha P. Chakrabarty for his guidance and invaluable comments on the manuscript.

\end{document}